\documentclass[12pt,oneside]{amsart}
\usepackage{amssymb, amsmath, amsthm}
\usepackage{epsfig}
\usepackage{epstopdf}
\usepackage{graphicx}
\usepackage{color}

\theoremstyle{plain}
\newtheorem{thm}{Theorem}

\newtheorem{prop}[thm]{Proposition}

\newtheorem{lemma}[thm]{Lemma}

\theoremstyle{definition}

      \makeatletter
      \def\@setcopyright{}
      \def\serieslogo@{}
      \makeatother

\begin{document}


   \title[Alternating Augmentations of Links]{Alternating Augmentations of Links}
   \author{Ryan C. Blair}

   \date{\today}


 \maketitle

\begin{abstract}
We show that one can interweave an unknot into any non-alternating
connected projection of a link so that the resulting augmented
projection is alternating.
\end{abstract}

In this paper, we will use the term link to mean a tame link
embedded in $S^3$. A link projection D is the image of a link
under a regular projection into $S^2$. D is a finite four-valent
graph in $S^2$.  If D is connected then each complement component
is a disk; we call the closure of one of these disks a region.
Define a labelling of the edges of D in the manner of Fig. 1.
Every edge of D receives two labels, one corresponding to each
end. An alternating edge is labelled with a plus and a minus while
a non-alternating edge receives two pluses or two minuses. Hence,
an alternating projection is one in which every edge is labelled
with a plus and a minus.  We call a non-alternating edge labelled
+ + a positive non-alternating edge and a non-alternating edge
labelled -- -- a negative non-alternating edge. Note that a
labelled knot projection is equivalent to the usual knot diagram
only with the crossing information stored as edge labels. In the
spirit of Adams\cite{a86}, we define an augmented link projection
of D to be the union of D with an unlink that projects to disjoint
simple closed curves in $S^2$.


We begin with a proposition originally outlined by
Thistlethwaite\cite{t88}.

\begin{figure}[h]
\centering \scalebox{.5}{\includegraphics{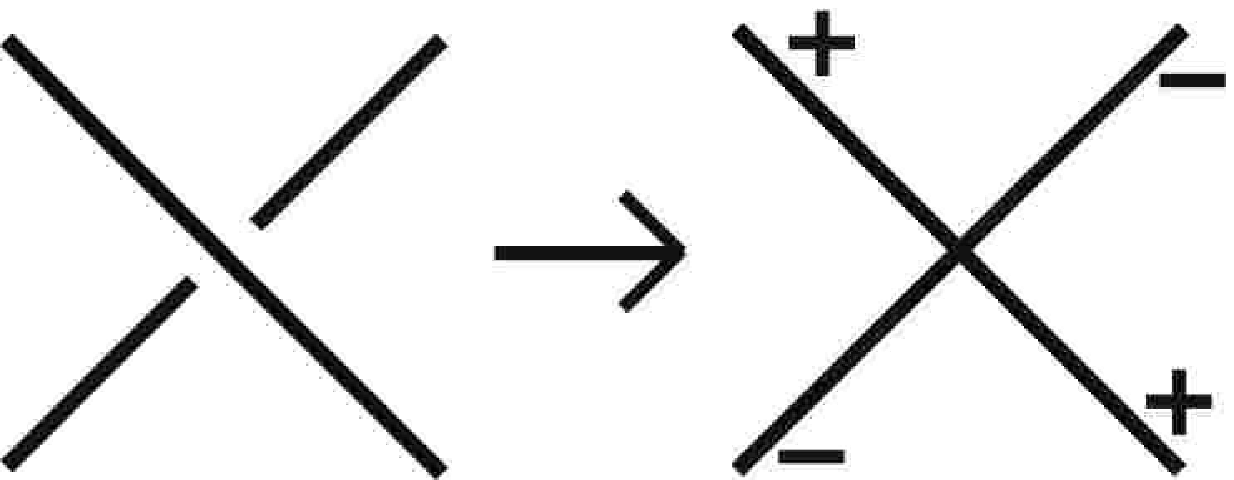}}
\caption{}\label{fig:labelling_small.eps}
\end{figure}

\begin{prop}
Any non-alternating connected link projection can be augmented so
that it becomes alternating.
\end{prop}
\begin{proof}
Choose a region R of the non-alternating projection D.  Let n be
the number of positive non-alternating edges of R.  Since each
vertex of $\partial$R contributes one + and one - label to the
edges of $\partial$R and each edge receives exactly two labels,
then $\partial$R contains exactly n negative non-alternating
edges. Let $\alpha$ be a non-alternating edge of $\partial$R.
After choosing an orientation of $\partial$R and again counting
the labels contributed by the vertices, the sign of the next
non-alternating edge in the direction of the orientation must be
opposite that of $\alpha$. The slogan here is that the sign of the
non-alternating edges of $\partial$R alternate.

With this knowledge, we can construct an augmented alternating
projection G from D as follows. Introduce vertices into the
interiors of all non-alternating edges of D. In every
non-alternating region R, we join pairs of such points together by
edges which are disjoint, lie in the interior of R and join
consecutive non-alternating edges as depicted in Fig. 2. We call
these new edges augmenting edges. Since the augmenting edges
connect non-alternating edges of opposite sign then there is a
consistent alternating labelling of the edges of G. In particular,
the end of an augmenting edge which bisects a positive
non-alternating edge of D receives a + label; the end that bisects
a negative non-alternating edge receives a -- label. Since the
augmenting edges never cross themselves, the closure of G-D
($cl(G-D)$) is a disjoint collection of simple closed curves
embedded in $S^2$. Because this process converts every
non-alternating edge of D into two alternating edges of G by
interweaving alternating unknots into D, G is an alternating
augmented projection of D.
\end{proof}

\begin{figure}[h]
\centering \scalebox{.5}{\includegraphics{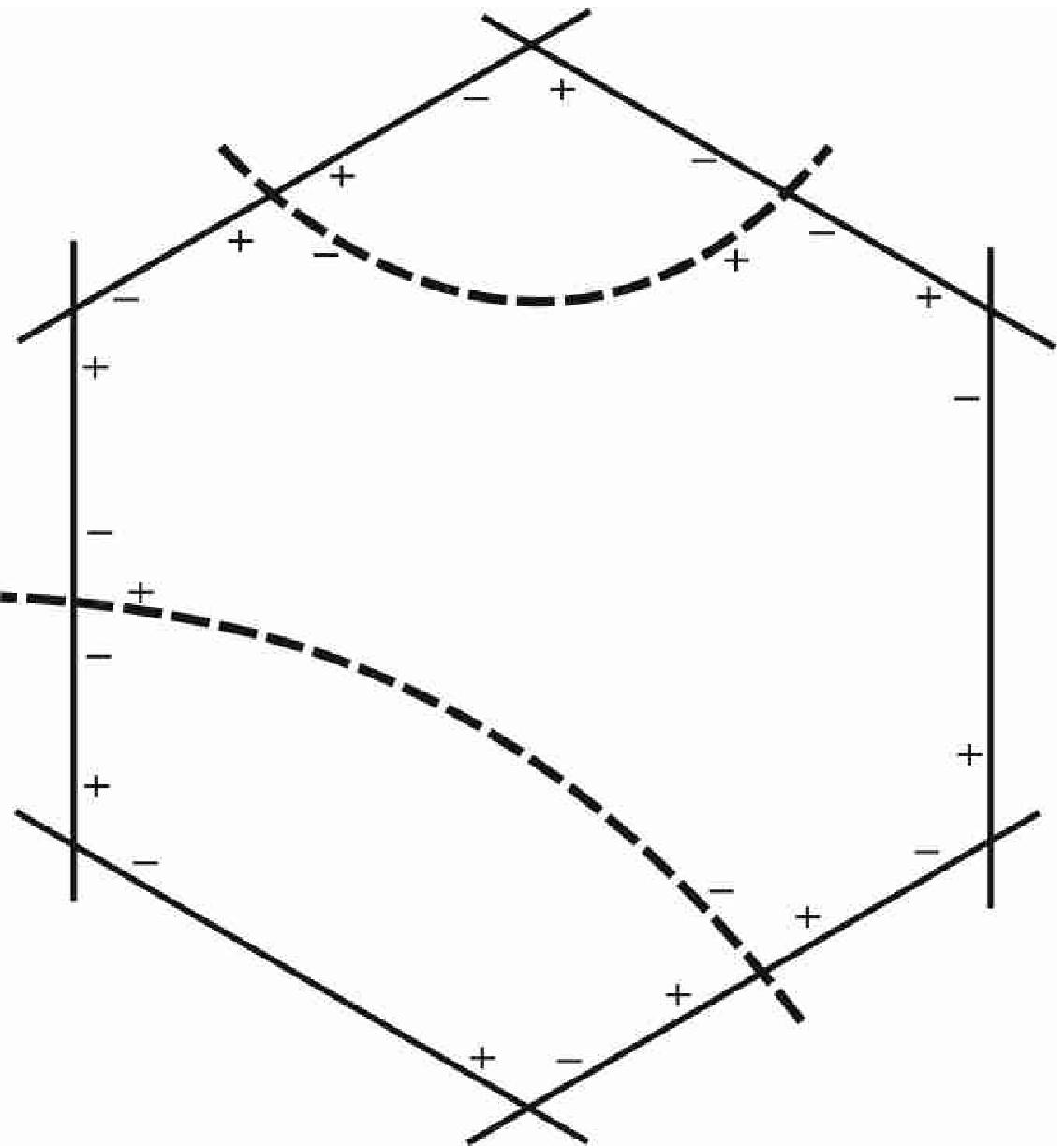}}
\caption{}\label{fig:moveI1.eps}
\end{figure}

Figure 3 illustrates an operation on link projections we call a
Type I move.

\begin{lemma}
Given an alternating connected projection of a link, a Type I move
results in an alternating link projection. \end{lemma}
\begin{figure}[h]
\centering \scalebox{.5} {\includegraphics{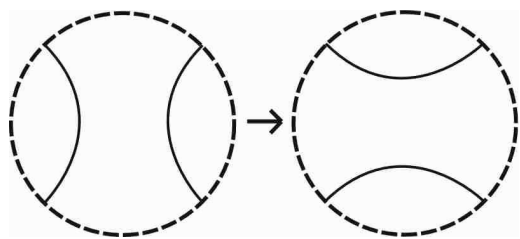}}
\caption{Type I Move}\label{fig:type1.eps}
\end{figure}
\begin{proof}
Let D be the link projection and $\alpha$ and $\beta$ be the
edges involved in the Type I move. $\alpha$ and $\beta$ are
boundary edges of a common region R. Since D is alternating, all
edges of R are labelled with both a plus and a minus.  Given two
consecutive edges of $\partial$R their shared vertex contributes a
plus label to one edge and a minus label to the other.  Thus, a
choice of label for a single edge determines the label of all the
edges of $\partial$R. In this way, an alternating label for
$\alpha$ determines the label for $\beta$, giving rise to the
following two possibilities.
\begin{figure}[h]
\centering \scalebox{.4} {\includegraphics{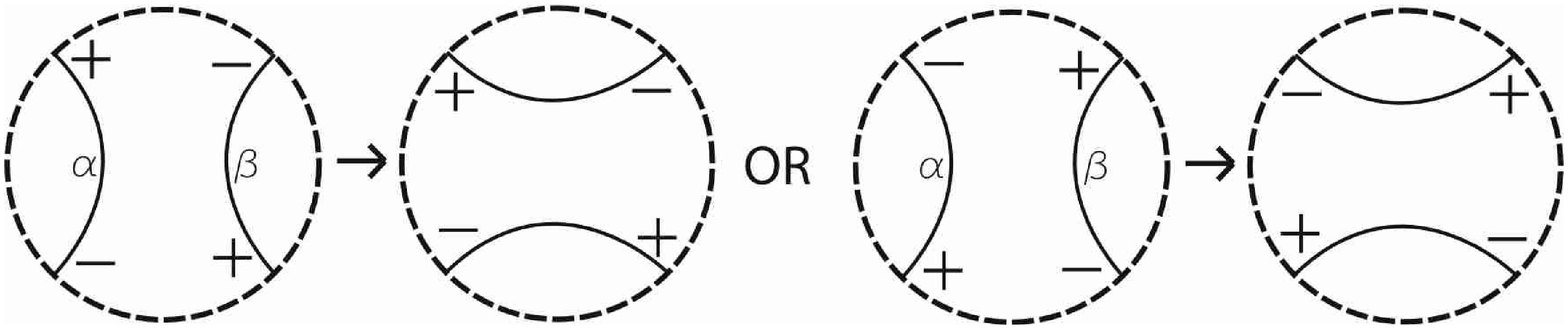}}
\caption{}\label{fig:type1marc.eps}
\end{figure}

In each case, the type I move preserves alternation.
\end{proof}

Figure 5 illustrates an operation on link projections we call a
Type II move.

\begin{lemma}
Given an alternating connected projection of a link, a Type II
move (after choosing labels incident to the new vertices) results
in an alternating link projection.
\end{lemma}
\begin{figure}[h]
\centering \scalebox{.5} {\includegraphics{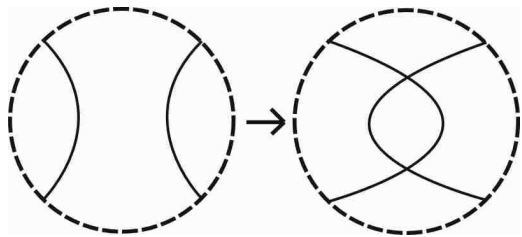}}
\caption{Type II Move}\label{fig:type2.eps}
\end{figure}
\begin{proof}
Let D be the link projection and $\alpha$ and $\beta$ be the edges
involved in the Type II move. We again use the fact that $\alpha$
and $\beta$ are boundary edges to a common region in D to deduce
that an alternating label for $\alpha$ determines the label for
$\beta$.  Hence, we have only the two following possibilities for
the labels of $\alpha$ and $\beta$.

\begin{figure}[h]
\centering \scalebox{.4} {\includegraphics{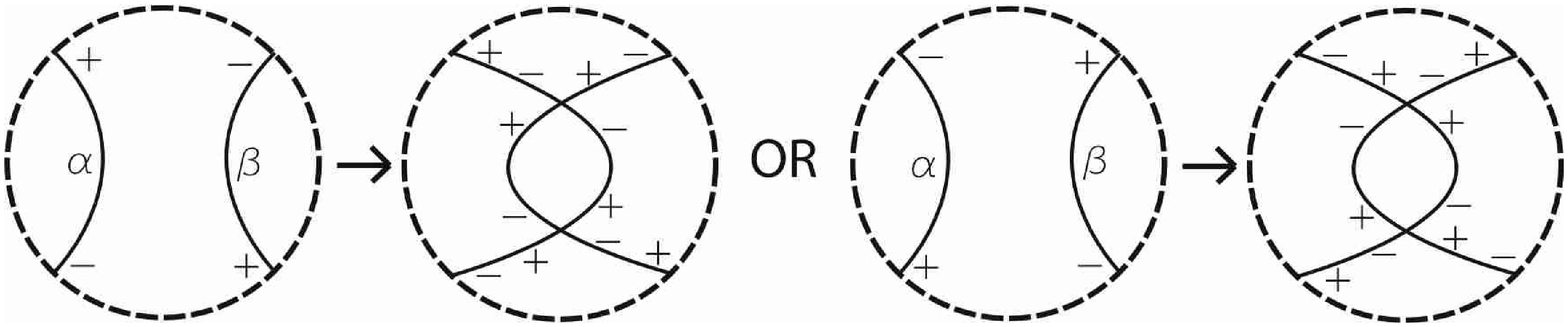}}
\caption{}\label{fig:type2marc.eps}
\end{figure}

In each case, as shown in Figure 6, we may choose a labelling of
the ends of the edges incident to the new vertices so that the
resulting diagrams are alternating.
\end{proof}

\begin{thm}
Given any connected projection of a non-alternating link, we can
augment the projection by adding a single unknotted component so
that the resulting link projection is alternating.
\end{thm}

\begin{proof}
Let D be a regular projection of a non-alternating link.  Create
G, an alternating augmented projection of D, as described in
Prop.1. Let $cl(G-D)= \bigcup_{1 \leq i \leq n}C_{i}$ where each
$C_{i}$ is a simple closed curve in $S^2$. If $n=1$ then there is
nothing to prove. If $n \geq 2$ then there is a path component A
of $S^2 - cl(G-D)$ whose closure has at least two distinct
boundary components $C_{i}$ and $C_{j}$.

If $C_{i}$ and $C_{j}$ contain boundary edges of a common region
in G then we may use a type I move to join $C_{i}$ and $C_{j}$
into a single simple closed curve in $S^2$. By Lemma 2 the
resulting projection is alternating.

\begin{figure}[h]
\centering \scalebox{.45}{\includegraphics{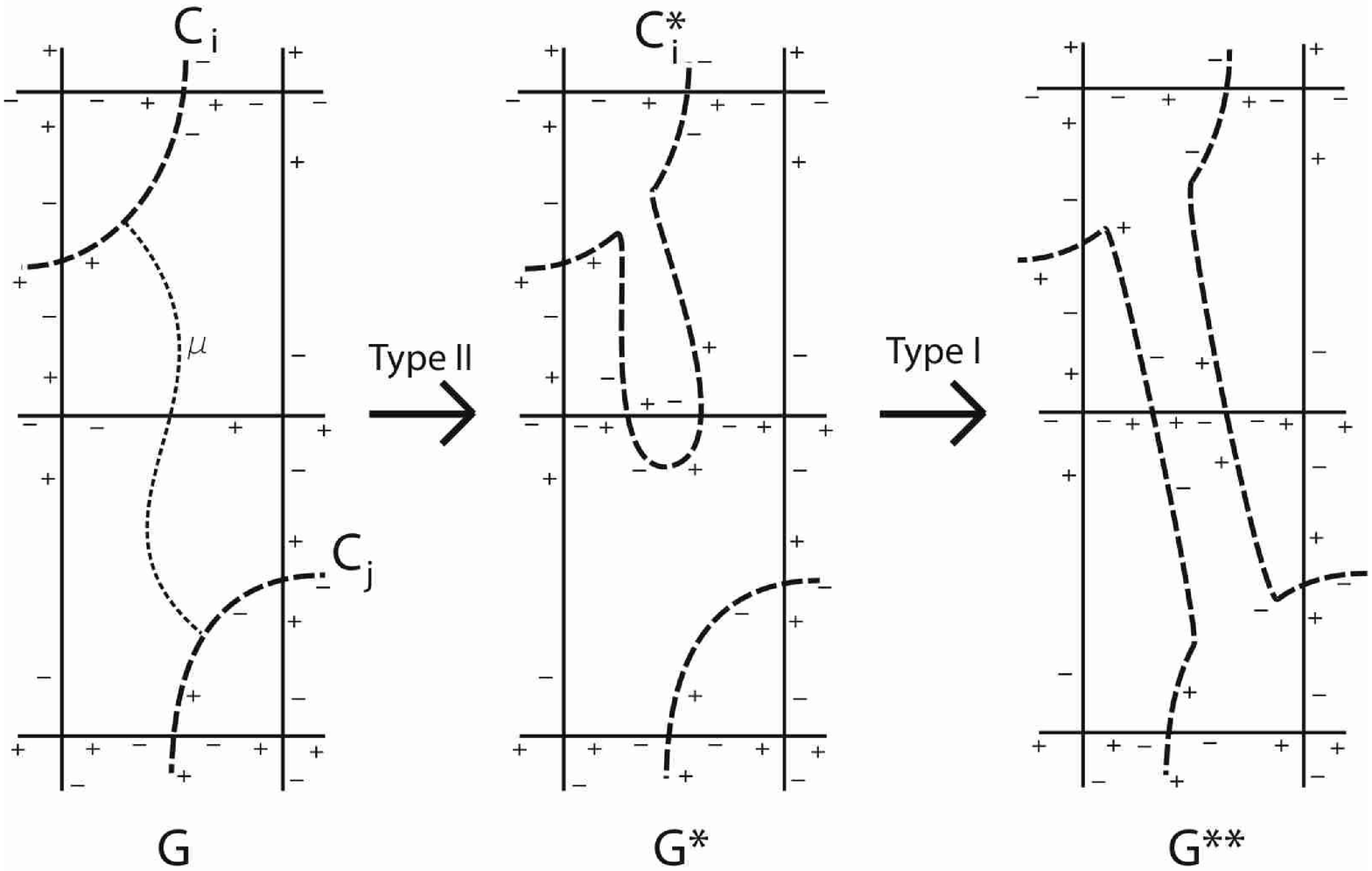}}
\caption{}\label{fig:movedia.eps}
\end{figure}

If $C_{i}$ and $C_{j}$ do not contain boundary edges of a common
region in G, then consider a path $\mu$ in A transverse to G so
that $\partial \mu = \{a,b\}$ for $a \epsilon C_{i}$ and $b
\epsilon C_{j}$.  We can propagate $C_{i}$ along $\mu$ using type
II moves as depicted in Fig. 7 until $C_{i}^{*}$(the image of
$C_{i}$ under type II moves) and $C_{j}$ contain boundary
components of a common region R. Call this projection $G^*$. Since
a type II move is an isotopy of $C_{i}$ in $S^2$ and $\mu$ was
restricted to A, $C_{i}^{*}$ is a simple closed curve which does
not intersect the other $C_{j}$. Hence, $G^{*}$ is an augmented
link of D. By Lemma 3, G is alternating implies we may choose
labels at the new vertices so that $G^*$ is alternating. We then
use the type I move to connect sum the disjoint simple closed
curves $C_{i}^{*}$ and $C_{j}$ into a single simple closed curve.
Call the resulting projection $G^{**}$. Since $G^{*}$ is
alternating so is $G^{**}$, by Lemma 2. Hence, $G^{**}$ is an
alternating augmented link of D with one less augmenting component
than G. Repeat this process until there is an alternating
augmented projection of D with exactly one augmenting component,
proving the theorem.
\end{proof}



\begin{thebibliography}{1}
\bibitem{a86} C. Adams, "Augmented alternating link complements are
hyperbolic," \emph{London Mathematical Society Lecture Notes
Series, 112: Low Dimensional Topology and Kleinian Groups},
pp.115-130, Cambridge University Press, Cambridge, 1986.
\bibitem{t88} Thistlethwaite, Morwen B, "An upper bound for the
breadth of the Jones polynomial," \emph{Math. Proc. Cambridge
Philos. Soc.} 103 (1988), no. 3, 451--456.
\end{thebibliography}
\end{document}